\documentclass[reqno,11pt]{amsart}
\usepackage{amsmath,amsxtra,amssymb,latexsym, amscd,amsthm}
\usepackage[unicode]{hyperref}
\usepackage{array,tabularx,longtable,multicol,indentfirst,fancybox,color}
\usepackage{graphicx}
\usepackage{multicol}
\usepackage{mathrsfs}

\advance\hoffset-1truecm\relax

\newtheorem{theorem}{Theorem}[section]
\newtheorem{definition}{Definition}
\newtheorem{lemma}{Lemma}[section]

\newtheorem{remark}{Remark}

\newtheorem{example}{Example}[section]
\numberwithin{equation}{section}

\def\R{\mathbb{R}}
\def\H{\mathscr{H}}
\def\B{\mathscr{B}}
\def\J{\mathcal{J}}
\def\T{\mathscr{T}}
\def\C{\mathcal{C}}

\begin{document}
	\title[Upper bound coefficient for convolution structure associated to Hartley--Bessel transform]{Upper bound coefficient for convolution structure \\associated to Hartley--Bessel transform}
	\date{\today.  Accepted by \textsc{Integral Transforms and Special Functions}}	
	\author[Trinh Tuan]{Trinh Tuan}
	\date{19 May 2026, 
	Accepted by \textsc{Integral Transforms and Special Functions}} 
\thanks{Published online: 25 May 2026.  \url{https://doi.org/10.1080/10652469.2026.2678565}}

\maketitle	
\begin{center}
\scriptsize	
	Department of Mathematics, Faculty of Natural Sciences, Electric Power University,\\ 235-Hoang Quoc Viet Rd.,  Hanoi, Vietnam.\\
	E-mail: \texttt{tuantrinhpsac@yahoo.com}
\end{center}	
\begin{abstract}
This paper is devoted to the study of a convolution structure denoted by $*_{\alpha}$, which is defined via the Hartley--Bessel transform. This concept was introduced in a recent work by F. Bouzeffour [\emph{J. Pseudo-Differ. Oper. Appl.}, 2024;15, Article 42]. We establish an analog of the Hausdorff--Young inequality for the Hartley--Bessel transform and convolution operator \( *_{\alpha} \). This leads to the convolution \( *_{\alpha} \) being uniformly bounded on the dual space. Moreover, in some special cases, our results yield a better upper-bound coefficient for the convolution \( *_{\alpha} \) than those previously obtained by Bouzeffour's result in [Theorem 4.4, \emph{J. Pseudo-Differ. Oper. Appl.}, 2024;15, Article 42].
Finally, we apply the convolution structure \( *_{\alpha} \) to study the solvability of a particular class of integral equations and provide a priori estimates for solutions under appropriate conditions.

\vskip 0.3cm

	\end{abstract}

\section{Introduction}
\subsection{Background and notions}
The Hartley transform, introduced in 1942 as an alternative to the classical Fourier transform, constitutes a symmetric integral representation originally developed for analyzing steady-state and transient phenomena in telephone transmission systems~\cite{Poularikas1996}. A distinguishing feature of the Hartley transform is its self-inverse nature; namely, the forward and inverse transforms are identical. Moreover, when applied to a real-valued function, the Hartley transform yields a real-valued frequency domain representation. This property confers a notable advantage over the classical Fourier transform, which typically involves complex-valued outputs even for real-valued inputs, thereby allowing the Hartley transform to operate entirely within the realm of real analysis \cite{Bracewell86}. Consequently, this transform has extensive applications across various fields of physics, engineering, and harmonic analysis \cite{TuanVK92,Hai92,yaku14,TuanMMA,tuan22,tuan2025refined}. 
Following \cite{Bracewell86}, the Hartley transform of a function \( f \) is defined as
\begin{equation}\label{eq1.1}
(\mathcal{H} f)(y) = \frac{1}{\sqrt{2\pi}} \int_{-\infty}^\infty f(x)\, \operatorname{cas}(xy)\ dx, \end{equation}
where \( y \in \R \) often represents an angular frequency in applications, and
$
\operatorname{cas}(x)= \cos(x) + \sin(x)
$
is known as basis kernel Hartley transform. By Euler's formula, this kernel can be expressed as
\begin{equation}\label{eq1.2}
\operatorname{cas}(x) = \mathrm{Re}\left((1 - i)e^{ix}\right) = \sum_{n=0}^\infty \frac{(-1)^{\frac{n(n-1)}{2}}}{n!} x^n.
\end{equation}
The function \(\operatorname{cas}\), defined in \eqref{eq1.2}, can be regarded as a real-analytic generalization of the exponential function. This interpretation becomes apparent when observing that \(\operatorname{cas}(x)\) is the unique \(C^\infty\)-solution to the  differential-reflection problem: $R \partial_x u(x) = \lambda u(x),$
 with initial condition $u(0) = 0$, see \cite{Bou21}.
Here \(\partial_x\) denotes the first-order derivative, and \(R\) is the reflection operator defined on functions \(f(x)\) by
$
(R f)(x) = f(-x).
$
An outstanding result established in~\cite{Bouze2014} asserts the following: If we consider the first-order difference-differential operator
$
\mathscr{L}_\alpha := R\left(\partial_x + \frac{\alpha}{x} \right) + \frac{\alpha}{x}$ with $\alpha > 0,
$
where \(\partial_x\) denotes differentiation with respect to the variable \(x\), and \(R\) is the reflection operator defined by \((Rf)(x) = f(-x)\). For any \(\lambda \in \mathbb{C}\), the eigenvalue problem
$
\mathscr{L}_\alpha y(x) = \lambda y(x)$ and $y(0) = 1,
$
admits a unique analytic solution 
\begin{equation}\label{eq1.3}
\mathcal{J}_\lambda (x, \alpha) =\B_{\alpha-\frac{1}{2}}(\lambda x)+\frac{\lambda x}{2\alpha+1}\B_{\alpha+\frac{1}{2}}(\lambda x),
\end{equation}
where $\B_{\alpha}(x)$  is the modified Bessel Function of the first kind of index $\alpha$ defined by \cite{Watson}
$$
\B_{\alpha}(x):=\sum_{n=0}^\infty \frac{(-1)^n}{n!(\alpha+1)_n}\left(\frac{x}{2}\right)^{2n}\ \text{with}\ \alpha>-1.
$$
Bessel functions of the first kind are solutions to the Bessel's differential equation. For non-negative integers or real parameters \( \alpha > 0 \), these functions remain finite at the origin \( x = 0 \). In contrast, when \( \alpha < 0 \) is non-integer, Bessel functions of the first kind diverge as \( x \to 0 \), (see Chapter 9 in \cite{Abra65}). 
Consequently, the function $\J_{\lambda}(x,\alpha)$ may be interpreted as a one-parameter extension of the classical \(\operatorname{cas}\) function. To emphasize this generalization, it is referred to as the Hartley--Bessel function (see \cite{Bloom95}). Subsequently in \cite{Bouze2024}, Bouzeffour  proposed the \textit{Hartley--Bessel transform} by replacing the
\(\operatorname{cas}\) kernel in~\eqref{eq1.1} with a  generalized Hartley--Bessel function. We recall some background on the Hartley--Bessel transform and its convolution structure, which has been recently introduced in \cite{Bouze2024}. Within the framework of this paper, we shall make frequent use of the weighted Lebesgue spaces, denoted by
$L^p_\alpha(\mathbb{R},\mu_\alpha(dx))$ be the space of measurable functions $f$ over $\mathbb{R}$ such that  $\int_{-\infty}^\infty |f(x)|^p \mu_\alpha(dx)$ is finite. For $1\leq p<\infty$ equipped with the
norm $$\|f\|_{L^p_\alpha(\mathbb{R},\mu_\alpha(dx))}=\bigg\{\int_{-\infty}^\infty |f(x)|^p \mu_\alpha(dx)\bigg\}^{1/p}.$$ In the case $p=\infty$, equipped with 
 $\|f\|_{L^\infty_\alpha(\mathbb{R},\mu_\alpha(dx))}= \underset{x\in\mathbb{R}}{\mathrm{ess\,sup}} |f(x)|.$ Here 
\( \mu_\alpha(dx) \) is a weighted measure given by
\begin{equation}\label{eq1.5}
\mu_\alpha(dx) = \frac{1}{2^{\alpha + \frac{1}{2}}\, \Gamma\left(\alpha + \frac{1}{2}\right)}\, |x|^{2\alpha}\, dx.
\end{equation}

\begin{definition}[$\mathscr{H}_\alpha$-transform]\label{dn1}
	Let \( \alpha>0 \), by \cite{Bouze2024} the Hartley--Bessel transform of a function \( f \in L^1_\alpha(\mathbb{R}, \mu_\alpha(dx)) \) is defined by
	\begin{equation}\label{eq1.6}
	\mathscr{H}_\alpha[f](\lambda) = \int_{-\infty}^\infty \mathcal{J}_\lambda (x, \alpha)\, f(x)\, \mu_\alpha(dx), \quad \lambda \in \mathbb{R},
	\end{equation}
	where \( \mathcal{J}_\lambda (x, \alpha) \) denotes the Hartley--Bessel function defined in~\eqref{eq1.3}, and \( \mu_\alpha(dx) \) determined by \eqref{eq1.5}.
\end{definition}
\noindent Moreover, by Proposition 6.4 in~\cite{Bouze2014} (refer also~\cite{Bouze2024}), the Hartley--Bessel transform admits an inversion formula. Specifically, if \( f \in L^1_\alpha(\mathbb{R}, \mu_\alpha(dx)) \) and \( \mathscr{H}_\alpha[f] \in L^1_\alpha(\mathbb{R}, \mu_\alpha(dx)) \), then :
\begin{equation}\label{biendoinguoc}
f(x) = \mathscr{H}_\alpha\left[ \mathscr{H}_\alpha[f](\lambda) \right](x), \ \text{a.e.\ on}\ \mathbb{R}. 
\end{equation}

\begin{definition}[Generalized translation operator $\T^x_\alpha$]	
	For any function \( f \in L^1_{\mathrm{loc}}(\mathbb{R}, \mu_\alpha(dx)) \), the generalized translation operator related to Hartley--Bessel operator is defined by
	\begin{equation}\label{eq1.7a}
	\mathscr{T}^x_\alpha f(t) = \int_{-\infty}^\infty f(s)\, \gamma_{x,t}^\alpha(ds),
	\end{equation}
	where \( \gamma_{x,t}^\alpha(ds) \) is a measure explicitly given by Definition 3.2 in \cite{Bouze2024}.
And \( L^1_{\mathrm{loc}}(\mathbb{R}, \mu_\alpha(dx)) \) denotes the space of Lebesgue measurable functions that are locally integrable on $\R$ with respect to the measure \( \mu_\alpha(dx) \).
\end{definition}
\noindent 		This translation operator $ \T^x_\alpha$ induces a new convolution structure, defined for all suitable functions \( f \) and \( g \), as follows:

\begin{definition}[Convolution structure $*_{\alpha}$]\label{dn2}
	The convolution product of two functions \( f \) and \( g \), associated with the Hartley--Bessel transform, is denoted as \( (f *_\alpha g) \) and defined by
	\begin{equation}\label{eq1.7}
	(f *_\alpha g)(x) = \int_{-\infty}^\infty \mathscr{T}^x_\alpha f(t)\, g(t)\, \mu_\alpha(dt), \quad x \in \mathbb{R},
	\end{equation}
	where \( \mathscr{T}^x_\alpha \) is the  translation operator defined by \eqref{eq1.7a}.
\end{definition}
\noindent
A notable result of Bouzeffour is Theorem 4.4 in \cite{Bouze2024},  provides boundedness properties of the convolution operator $(f *_\alpha g)$, asserts that:
\textbf{(i)} Let \( \alpha > 0 \). For \( 1 \leq p, q, r \leq \infty \) satisfy the relation
	$
	\frac{1}{p} + \frac{1}{q} - 1 = \frac{1}{r}.
	$
	For all functions $f\in L_\alpha^p(\mathbb{R},\mu_\alpha(dx))$, $g\in L_\alpha^q(\mathbb{R},\mu_\alpha(dx))$, the convolution product $( f *_\alpha g )$ is well-defined in \( L_\alpha^r(\mathbb{R},\mu_\alpha(dx)) \), and the following inequality holds
	\begin{equation}\label{eq1.8}
	\left\|f \ast_\alpha g\right\|_{L_\alpha^r (\mathbb{R},\mu_\alpha(dx))} \leq 4\|f\|_{L_\alpha^p(\mathbb{R},\mu_\alpha(dx))} \|g\|_{L_\alpha^q(\mathbb{R},\mu_\alpha(dx))}.
	\end{equation}
	\textbf{(ii)} Let \( \alpha > 0 \).  For \( 1 \leq p, q, r \leq 2 \) satisfy condition 
	$
	\frac{1}{p} + \frac{1}{q} - 1 = \frac{1}{r}.
	$
	Then for all $f\in L_\alpha^p(\mathbb{R},\mu_\alpha(dx))$, $g\in L_\alpha^q(\mathbb{R},\mu_\alpha(dx))$, the Hartley--Bessel transform satisfies factorization identity 
	\begin{equation}\label{eq1.9a}
	\mathscr{H}_\alpha[f *_\alpha g]= \mathscr{H}_\alpha [f]	 \mathscr{H}_\alpha [g].
	\end{equation}
In addition, the convolution $\ast_\alpha$ is associative in \( L^1_{\alpha} (\mathbb{R}, \mu_\alpha(dx)) \).	
Notice that, for the case $p=q \in [1,2]$, by applying the $\mathscr{H}_\alpha$-transform on both sides of the expression \eqref{eq1.9a} and combining \eqref{biendoinguoc}, we infer that 
\begin{equation}\label{eq1.9}
f\ast_\alpha g= \mathscr{H}_\alpha\Big(\mathscr{H}_\alpha[f] \mathscr{H}_\alpha[g]\Big)\ \text{with}\ \alpha>0.
\end{equation}
\begin{remark}\label{daisobanach}
\textup{It is worth emphasizing that, in light of asserts (ii) of Theorem 4.4 in \cite{Bouze2024}, if $p=q=r=1$, for any function $f,g \in L^1_\alpha(\mathbb{R}, \mu_\alpha(dx))$ then $f *_\alpha g \in L^1_\alpha(\mathbb{R}, \mu_\alpha(dx))$ is well-defined and the following estimate always occurs
$$
	\left\|f *_\alpha g\right\|_{L^1_\alpha(\mathbb{R}, \mu_\alpha(dx))} \leq \|2f\|_{L^1_\alpha(\mathbb{R}, \mu_\alpha(dx))} \|2g\|_{L^1_\alpha(\mathbb{R}, \mu_\alpha(dx))}.
$$
	This implies that the space \( L^1_\alpha(\mathbb{R}, \mu_\alpha(dx)) \), when equipped with the multiplication, is convolution as defined by  \eqref{eq1.7}, forms a Banach algebra structure \cite{1964,Naimark}; (also refer \cite{Dales00}).}
\end{remark}

\noindent

\subsection{Purpose and organization}
Our main result of this paper is to get another Hausdorff--Young inequality similar to \eqref{eq1.8} but in the space $L^{r_1}_\alpha (\mu_\alpha(dx))$ where $1/ r +1/r_1 =1$, instead of $L^r_\alpha(\mathbb{R}, \mu_\alpha(dx))$. In the special case when $r=r_1 =2$, our inequality reduces to \eqref{eq1.8}, but with a better upper bound constant $\sqrt 2$ (instead of $4$ as  Bouzeffour's result, Theorem 4.4 in \cite{Bouze2024}).

This paper is organized into three sections. As stated in the Abstract, the principal results of this paper are established in Section~\ref{sec2}. The final part, Section~\ref{sec3}, is devoted to exploring the applicability of the convolution structure~\eqref{eq1.7} to the solvability of a specific class of integral equations. This analysis relies on the theoretical framework developed in the preceding sections, together with the Wiener--Lévy invertibility criterion in \cite{1964,Naimark} and related convolution inequalities.
\section{Statement of main result}\label{sec2}
\noindent For real numbers $r, r_1>1$ , they are called conjugate indices (or Hölder conjugates) if
$1/r +1/r_1=1.
$
Formally, we also define $r_1=\infty$ as conjugate to $r=1$ and vice versa.
Conjugate indices are used in Hölder's inequality, as well as Young's inequality for convolution products. If $r, r_1>1$ are conjugate exponents, then the Lebesgue spaces  $L^r$ and $L^{r_1}$ are dual to each other in the sense of functional analysis. By Theorem~4.4 in~\cite{Bouze2024}, it can be asserted that: Let \( \alpha > 0\) and \( p, q, r \in [1, \infty] \) satisfy the relation  
$
1/p+1/q=1 +1/r,
$ 
then the convolution structure defined in~\eqref{eq1.7} yields a bilinear mapping  
$$
\ast_\alpha:L^p_\alpha(\mathbb{R}, \mu_\alpha(dx)) \times L^q_\alpha(\mathbb{R}, \mu_\alpha(dx)) \to L^r_\alpha(\mathbb{R}, \mu_\alpha(dx)),
$$ 
such that the map \( (f, g) \mapsto f \ast_\alpha g \) is continuous in each argument, for \( f \in L^p_\alpha(\mathbb{R}, \mu_\alpha(dx)) \) and \( g \in L^q_\alpha(\mathbb{R}, \mu_\alpha(dx)) \). Furthermore, the convolution operator \eqref{eq1.7}
is uniformly bounded on $L^r_\alpha(\mathbb{R}, \mu_\alpha(dx))$ space. 
The question that arises here is whether the convolution  \eqref{eq1.7} remains uniformly bounded when acting on the dual space of \( L^r_\alpha(\mathbb{R}, \mu_\alpha(dx)) \)? To answer this question, we first require the following lemma.
\begin{lemma}[Hausdorff--Young inequality for $\mathscr{H}_\alpha$-transform]\label{lemma21}
Let $p \in [1,2]$ satisfy $1/p + 1/p_1 =1$. For any functions $f\in L_\alpha^p(\mathbb{R},\mu_\alpha(dx))$, the following estimation holds
\begin{equation}\label{eq2.1}
\left\|\mathscr{H}_\alpha[f]\right\|_{L_\alpha^{p_1}(\mathbb{R},\mu_\alpha(dx))} \leq (\sqrt{2})^{\frac{2}{p} -1} \|f\|_{L_\alpha^p(\mathbb{R},\mu_\alpha(dx))}. 
\end{equation}
\end{lemma}
\begin{proof}
	Based on Corollary~2.3 in~\cite{Bouze2024}, for all \( \alpha > 0 \) and \( \lambda, x \in \mathbb{R} \), the Hartley--Bessel kernel satisfies the uniform bound
	$
	|\J_\lambda(x, \alpha)| \leq \sqrt{2}.
	$
	This pointwise estimate plays a crucial role in analyzing the properties of the convolution structure. Indeed, from the definition of convolution ~\eqref{eq1.7}, $\forall f\in L_\alpha^1 (\mathbb{R},\mu_\alpha(dx))$ we infer that
	\begin{align*}
	|\mathscr{H}_\alpha[f](\lambda)|&\leq \int_{-\infty}^\infty |\J_\lambda(x,\alpha)| |f(x)| \mu_\alpha (dx)\\
	&\leq \sqrt{2}\int_{-\infty}^\infty |f(x)| \mu_\alpha (dx)
	= \sqrt{2}\|f\|_{L_\alpha^1(\mathbb{R},\mu_\alpha(dx))}<\infty.
	\end{align*}
This means that \begin{equation}\label{eq2.2}
\|\mathscr{H}_\alpha[f]\|_{L_\alpha^\infty(\mathbb{R},\mu_\alpha(dx))} \leq \sqrt{2}\|f\|_{L_\alpha^1(\mathbb{R},\mu_\alpha(dx))},
\end{equation}
where $\|\varphi\|_{L_\alpha^\infty(\mathbb{R},\mu_\alpha(dx))} = \underset{x\in\mathbb{R}}{\mathrm{ess\,sup}}\,|\varphi(x)|$. Therefore, the $\mathscr{H}_\alpha$-transform \eqref{eq1.6} is the bounded operator from $L_\alpha^1(\mathbb{R},\mu_\alpha(dx)) \to L_\alpha^\infty(\mathbb{R},\mu_\alpha(dx)).$  Moreover, there exists a unique isometric extension of the operator \( \mathscr{H}_\alpha \) to the whole space \( L_\alpha^2(\mathbb{R},\mu_\alpha(dx)) \), which coincides with \( \mathscr{H}_\alpha \) on the dense subspace $L_\alpha^1(\mathbb{R},\mu_\alpha(dx))\cap L_\alpha^2(\mathbb{R},\mu_\alpha(dx))$, refer \cite{Bouze2014} (also see \cite{Bouze2024}, Statement (3) in Theorem 4.1). This leads to for any functions $f\in L_\alpha^1(\mathbb{R},\mu_\alpha(dx)) \cap L_\alpha^2(\mathbb{R},\mu_\alpha(dx))$ then $\mathscr{H}_\alpha[f]\in L_\alpha^2(\mathbb{R},\mu_\alpha(dx))$ and we get
\begin{equation}\label{eq2.3}
 \|\mathscr{H}_\alpha[f]\|_{L_\alpha^2(\mathbb{R},\mu_\alpha(dx))} = \|f\|_{L_\alpha^2(\mathbb{R},\mu_\alpha(dx))}.
\end{equation} 
This implies that $\mathscr{H}_\alpha$ is bounded from $L_\alpha^2(\mathbb{R},\mu_\alpha(dx))$ to $ L_\alpha^2(\mathbb{R},\mu_\alpha(dx))$. Now, by applying Riesz-Thorin's interpolation (see \cite[Theorem 1.19]{Duoandikoetxea01} also refer \cite{Stein1972Weiss}, Chapter 5), for $\alpha >0$  we deduce $\mathscr{H}_\alpha$ is a bounded linear operator from $L_\alpha^p(\mathbb{R},\mu_\alpha(dx)) \to L_\alpha^{p_1}(\mathbb{R},\mu_\alpha(dx))$, and $$
\|\mathscr{H}_\alpha[f]\|_{L_\alpha^{p_1}(\mathbb{R},\mu_\alpha(dx))} \leq \left(\sqrt{2}\right)^{\theta}\|f\|_{L_\alpha^p(\mathbb{R},\mu_\alpha(dx))},
$$ where $0<\theta<1$, defined by $\frac{1}{p}=\frac{\theta}{1}+\frac{1-\theta}{2}$, $\forall p \in [1,2]$. It follows that $\theta=\frac{2}{p} -1$, and we obtain an estimation as in the conclusion 
of this lemma with an upper bound coefficient is $\left(\sqrt{2}\right)^{\frac{2}{p} -1}$.
\end{proof}
\begin{remark}
\textup{We examine the estimates at endpoints of the interval \( p \in [1,2] \) as established in Lemma~\ref{lemma21}. When \( p = 1 \), it follows that \( p_1 = \infty \), and the inequality~\eqref{eq2.1} reduces to~\eqref{eq2.2}. On the other hand, when \( p = 2 \), we have \( p_1 = 2 \), and inequality~\eqref{eq2.1} becomes to~\eqref{eq2.3}. Consequently, with \( p \in [1, 2] \) then \( \left(\sqrt{2}\right)^{\frac{2}{p} - 1} \) is  the sharp optimal in~\eqref{eq2.1} for all functions \( f \in L^p_\alpha(\mathbb{R}, \mu_\alpha(dx)) \).}
\end{remark}

\begin{theorem}[\textsc{Main Theorem}]\label{mainTHEOREM}
	Suppose that $p,q,r \in [1,2]$ are real numbers satisfying the relations $$1/p +1/p_1 = 1/q + 1/q_1=1/r +1/r_1 =1\ \  \text{and}\ \ 1/p_1 +1/q_1 = 1/r.
	$$ Then, for any functions  $f\in L_\alpha^p(\mathbb{R},\mu_\alpha(dx))$, $g\in L_\alpha^q(\mathbb{R},\mu_\alpha(dx))$, with $\alpha >0$, the following estimate holds:
	\begin{equation}\label{eq2.4}
	\|f \ast_\alpha g\|_{L_\alpha^{r_1}(\mathbb{R},\mu_\alpha(dx))} \leq \C_{p,q,r}\|f\|_{L_\alpha^p(\mathbb{R},\mu_\alpha(dx))} \|g\|_{L_\alpha^q(\mathbb{R},\mu_\alpha(dx))},
	\end{equation}	
	where the coefficient $\C_{p,q,r}$ given by $\C_{p,q,r}:=(\sqrt{2})^{\left(\frac{2}{p}+\frac{2}{q}+\frac{2}{r}-3\right)}$ is the upper bound on the right-hand side of the
	inequality \eqref{eq2.4} and depends explicitly on the numbers $p, q, r$.
\end{theorem}
\begin{proof}
	From the assumption  $1/p_1 +1/q_1 = 1/r$, we deduce $\frac{r}{p_1}+\frac{r}{q_1}=1$.  Applying the Hölder inequality we have
\begin{equation*}
\begin{aligned}
&\int_{-\infty}^\infty \left|\mathscr{H}_\alpha[f](\lambda) \mathscr{H}_\alpha[g](\lambda)\right|^r \mu_\alpha (d\lambda)\\&
=\int_{-\infty}^\infty \left(\big|\mathscr{H}_\alpha[f](\lambda) \mathscr{H}_\alpha[g](\lambda)\big|\mu_\alpha^{\frac{1}{r}}\right)^r(d\lambda)\\
&=\int_{\infty}^\infty\left(\left|\mathscr{H}_\alpha[f](\lambda)\right|^r\mu_\alpha^{\frac{r}{p_1}}\right)\left(\left|\mathscr{H}_\alpha[g](\lambda)\right|^r\mu_\alpha^{\frac{r}{q_1}}\right)(d\lambda),\\
&\leq \left\{\int_{-\infty}^\infty \left[\left(\left|\mathscr{H}_\alpha[f](\lambda)\right|\mu_\alpha^{\frac{1}{p_1}}\right)^r\right]^{\frac{p_1}{r}}(d\lambda)\right\}^{\frac{r}{p_1}}\times \left\{\int_{-\infty}^\infty \left[\left(\left|\mathscr{H}_\alpha[g](\lambda)\right|\mu_\alpha^{\frac{1}{q_1}}\right)^r\right]^{\frac{q_1}{r}} (d\lambda)\right\}^{\frac{r}{q_1}}\\
&=\left\{\int_{-\infty}^\infty\left|\mathscr{H}_\alpha[f](\lambda)\right|^{p_1}\mu_\alpha (d\lambda)\right\}^{\frac{r}{p_1}} \times \left\{\int_{-\infty}^\infty\left|\mathscr{H}_\alpha[g](\lambda)\right|^{q_1}\mu_\alpha(d\lambda)\right\}^{\frac{r}{q_1}}.
\end{aligned}
\end{equation*}
This yields
\begin{equation}\label{eq2.5}
\int_{-\infty}^\infty \left|\mathscr{H}_\alpha[f](\lambda) \mathscr{H}_\alpha[g](\lambda)\right|^r \mu_\alpha (d\lambda) \leq \left\|\mathscr{H}_\alpha[f]\right\|_{L_\alpha^{p_1}(\mathbb{R},\mu_\alpha(dx))}^r\left\|\mathscr{H}_\alpha[g]\right\|_{L_\alpha^{q_1}(\mathbb{R},\mu_\alpha(dx))}^{r}
\end{equation}
Applying the estimate \eqref{eq2.1} to inequality \eqref{eq2.5}, we obtain
\begin{equation}\label{eq2.6}
\begin{aligned}
&\int_{-\infty}^\infty\big|\mathscr{H}_\alpha[f](\lambda) \mathscr{H}_\alpha[g](\lambda)\big|^r \mu_\alpha (dx)
\leq\left\|\mathscr{H}_\alpha[f]\right\|_{L_\alpha^{p_1}(\mathbb{R},\mu_\alpha(dx))}^r \left\|\mathscr{H}_\alpha[g]\right\|_{L_\alpha^{q_1}(\mathbb{R},\mu_\alpha(dx))}^r\\
&\leq (\sqrt{2})^{r\left(\frac{2}{p}-1\right)}\|f\|_{L_\alpha^p(\mathbb{R},\mu_\alpha(dx))}^r\times (\sqrt{2})^{r\left(\frac{2}{q}-1\right)}\|g\|_{L_\alpha^q(\mathbb{R},\mu_\alpha(dx))}^r.
\end{aligned}
\end{equation}
This shows that, for any functions $f\in L_\alpha^p(\mathbb{R},\mu_\alpha(dx))$ and $g\in L_\alpha^q(\mathbb{R},\mu_\alpha(dx))$, then $$\mathscr{H}_\alpha[f](\lambda)\mathscr{H}_\alpha[g](x) \in L_\alpha^r(\mathbb{R},\mu_\alpha(dx)).$$ Coupling \eqref{eq1.9}, \eqref{eq2.1} and \eqref{eq2.6}, we infer that
\begin{align*}
\left\|f \ast_\alpha g\right\|_{L_\alpha^{r_1}(\mathbb{R},\mu_\alpha(dx))}^r&=\left\|\mathscr{H}_\alpha\left[\mathscr{H}_\alpha[f](\lambda) \mathscr{H}_\alpha[g](\lambda)\right]\right\|_{L_\alpha^{r_1}(\mathbb{R},\mu_\alpha(dx))}^r\\
&\leq (\sqrt{2})^{r\left(\frac{2}{r}-1\right)}\left\|\mathscr{H}_\alpha[f](\lambda) \mathscr{H}_\alpha[g](\lambda)\right\|_{L_\alpha^r(\mathbb{R},\mu_\alpha(dx))}^r\\
&=(\sqrt{2})^{r\left(\frac{2}{r}-1\right)}\int_{-\infty}^\infty \left|\mathscr{H}_\alpha[f](\lambda) \mathscr{H}_\alpha[g](\lambda)\right|^r \mu_\alpha (dx)\\
&\leq (\sqrt{2})^{r\left(\frac{2}{p}+\frac{2}{q}+\frac{2}{r}-3\right)}\|f\|_{L_\alpha^p(\mathbb{R},\mu_\alpha(dx))}^r \|g\|_{L_\alpha^q(\mathbb{R},\mu_\alpha(dx))}^r.
\end{align*}
This completes the proof.
\end{proof}
\noindent Theorem~\ref{mainTHEOREM} establishes that the convolution structure \( \ast_\alpha \), for \( \alpha > 0 \), defines a bounded bilinear operator on the space \( L_\alpha^{r_1}(\mathbb{R}, \mu_\alpha(dx)) \) for all functions \( f \in L_\alpha^p(\mathbb{R}, \mu_\alpha(dx)) \) and \( g \in L_\alpha^q(\mathbb{R}, \mu_\alpha(dx)) \). Moreover, under the conditions on the exponents specified in Theorem~\ref{mainTHEOREM}, the associated constant satisfies
\[ 
\mathcal{C}_{p,q,r} = (\sqrt{2})^{\left(\frac{2}{p}+\frac{2}{q}+\frac{2}{r}-3\right)} < 2\sqrt{2}.
\]
This yields the following inequality:
\begin{equation*}\label{eq2.7}
\left\|f \ast_\alpha g\right\|_{L_\alpha^{r_1}(\mathbb{R}, \mu_\alpha(dx))} 
< 2\sqrt{2} \, \|f\|_{L_\alpha^p(\mathbb{R}, \mu_\alpha(dx))} \, \|g\|_{L_\alpha^q(\mathbb{R}, \mu_\alpha(dx))},
\end{equation*}
which confirms that the convolution operator \eqref{eq1.7} is uniformly bounded in \( L_\alpha^{r_1}(\mathbb{R}, \mu_\alpha(dx)) \). This result provides an affirmative answer to the question raised at the beginning of this section. In what follows, we analyze several special cases of Theorem~\ref{mainTHEOREM}.

\begin{remark}\label{rem3}
\textup{\textbf{(A).} For the case \( r = 1 \), then necessarily \( r_1 = \infty \). Let us consider the choice \( p = p_1 = q = q_1 = 2 \), which clearly satisfies the index conditions stated in Theorem~\ref{mainTHEOREM}. In this case, the convolution inequality \eqref{eq2.4} becomes
	\begin{equation}\label{eq2.8}
	\left\|f \ast_\alpha g\right\|_{L_\alpha^\infty(\mathbb{R}, \mu_\alpha(dx))} \leq \sqrt{2} \, \|f\|_{L_\alpha^2(\mathbb{R}, \mu_\alpha(dx))} \, \|g\|_{L_\alpha^2(\mathbb{R}, \mu_\alpha(dx))},
	\end{equation}
	where the \( L^\infty \)-norm is defined as
	$
	\|\varphi\|_{L_\alpha^\infty(\mathbb{R}, \mu_\alpha(dx))} := \underset{x \in \mathbb{R}}{\mathrm{ess\,sup}} \, |\varphi(x)|.
$	
	Let us now compare this result with inequality~\eqref{eq1.8}, which corresponds to Item (i) of Theorem~4.4 in~\cite{Bouze2024}. When taking \( p = q = 2 \in [1, \infty] \), the index relation \( \frac{1}{p} + \frac{1}{q} = 1 + \frac{1}{r} \) implies \( r = \infty \), and Theorem~4.4 yields the following estimate:
	\[
	\left\|f \ast_\alpha g\right\|_{L_\alpha^\infty(\mathbb{R}, \mu_\alpha(dx))} \leq 4 \, \|f\|_{L_\alpha^2(\mathbb{R}, \mu_\alpha(dx))} \, \|g\|_{L_\alpha^2(\mathbb{R}, \mu_\alpha(dx))}.
	\]
	It is evident that the constant \( (\sqrt{2}) \) in the inequality \eqref{eq2.8} provides a better bound than the constant \( (4) \) appearing in Theorem~4.4 of~\cite{Bouze2024}. Therefore, in this specific case, our result improves upon the known estimate in the literature.\\
\textbf{(B).} 
When \( r = 2 \), we have \( r_1 = 2 \) by duality. Under the index conditions prescribed in Theorem~\ref{mainTHEOREM}, one may consider the following two admissible configurations of exponents:
\begin{itemize}
	\item \( p = p_1 = 2 \), \( q = 1 \), and \( q_1 = \infty \), or
	\item \( q = q_1 = 2 \), \( p = 1 \), and \( p_1 = \infty \).
\end{itemize}
	Applying inequality~\eqref{eq2.4} in each case yields the following norm estimates, respectively:
	\begin{align}
	\left\|f \ast_\alpha g\right\|_{L_\alpha^2(\mathbb{R}, \mu_\alpha(dx))} &\leq \sqrt{2} \, \|f\|_{L_\alpha^1(\mathbb{R}, \mu_\alpha(dx))} \, \|g\|_{L_\alpha^2(\mathbb{R}, \mu_\alpha(dx))}, \quad \text{and} \label{eq2.9}\\
	\left\|f \ast_\alpha g\right\|_{L_\alpha^2(\mathbb{R}, \mu_\alpha(dx))} &\leq \sqrt{2} \, \|f\|_{L_\alpha^2(\mathbb{R}, \mu_\alpha(dx))} \, \|g\|_{L_\alpha^1(\mathbb{R}, \mu_\alpha(dx))} \label{eq2.10}.
	\end{align}
	To assess the sharpness of these inequalities, we compare them with the corresponding estimates derived from Theorem~4.4 of~\cite{Bouze2024}. In the same setting where (\( r = q = 2 \); \( p = 1 \)), or (\( r = p = 2 \); \( q = 1 \)), the result from~\cite{Bouze2024} asserts the upper bounds:
	$$\begin{aligned}
		\left\|f \ast_\alpha g\right\|_{L_\alpha^2(\mathbb{R},\mu_\alpha(dx))} &\leq 4 \, \|f\|_{L_\alpha^1(\mathbb{R},\mu_\alpha(dx))} \, \|g\|_{L_\alpha^2(\mathbb{R},\mu_\alpha(dx))},\quad \text{and}\\
		\left\|f \ast_\alpha g\right\|_{L_\alpha^2(\mathbb{R},\mu_\alpha(dx))} &\leq 4 \, \|f\|_{L_\alpha^2(\mathbb{R},\mu_\alpha(dx))} \, \|g\|_{L_\alpha^1(\mathbb{R},\mu_\alpha(dx))}.
	\end{aligned}$$
Comparing these results, it is evident that the constants \( (\sqrt{2}) \) appearing in inequalities~\eqref{eq2.9} and~\eqref{eq2.10} provide a substantial improvement over the constant \( (4) \) in Theorem~4.4 of~\cite{Bouze2024}, thereby demonstrating the sharpness of our result in this context.}
\end{remark}

\section{Application and example}\label{sec3}
\noindent 
In this section, we investigate a class of integral equations associated with the convolution structure defined in \eqref{eq1.7}. Our objective is to derive sufficient conditions under which the existence and uniqueness of closed-form solutions can be guaranteed. Specifically, we consider the integral equation of the form
\begin{equation}\label{eq3.3}
f(x) + \int_{-\infty}^\infty K(\alpha, x, t) f(t) dt = \left(g \ast_\alpha h\right)(x)\  \text{with}\ \alpha >0,
\end{equation}
where \( f(x) \) is an unknown function to be determined, \( K(\alpha, x, t) \) is a given kernel, and \( g, h \in L_\alpha^1(\mathbb{R}, \mu_\alpha(dx)) \) are prescribed functions. We will show the conditions for solvability in $L_\alpha^1(\mathbb{R}, \mu_\alpha(dx))$ of equation \eqref{eq3.3} for case the kernel $K(\alpha,x,t):= \T_\alpha^x g(t) \mu_\alpha$ with $(\alpha>0)$, where $\T_\alpha^x$ is a translation operator defined by \eqref{eq1.6} and $\mu_\alpha$ is a weighted measure given by \eqref{eq1.5}. By definition of $\ast_\alpha$ in \eqref{eq1.7}, then \eqref{eq3.3} can be rewritten in convolution form 
\begin{equation}\label{eq3.4}
f(x)+ \left(f \ast_\alpha g\right)(x) = \left(g \ast_\alpha h\right)(x)\ \text{with}\ \alpha >0.
\end{equation}

\begin{theorem}\label{theo3.1}
Let \( g, h \in L_\alpha^1(\mathbb{R}, \mu_\alpha(dx)) \) are given functions. For the equation \eqref{eq3.4} to be solvable in the space \( L_\alpha^1(\mathbb{R}, \mu_\alpha(dx)) \), a sufficient condition is that
$
1 + \mathscr{H}_\alpha[g](\lambda) \neq 0,\  \forall \lambda \in \mathbb{R}.
$
Under this condition, the solution \( f \in L_\alpha^1(\mathbb{R}, \mu_\alpha(dx)) \) exists uniquely and is given almost everywhere on \( \mathbb{R} \) by formula
$$
f(x) = \left(\ell \ast_\alpha h\right)(x),
$$
where the function \( \ell \in L_\alpha^1(\mathbb{R}, \mu_\alpha(dx)) \) is defined through
$
\mathscr{H}_\alpha[\ell](\lambda) = \frac{\mathscr{H}_\alpha[g](\lambda)}{1 + \mathscr{H}_\alpha[g](\lambda)}, \forall \lambda \in \mathbb{R}.
$
Moreover, the solution \( f \) satisfies the estimate
$$
\|f\|_{L_\alpha^1(\mathbb{R}, \mu_\alpha(dx))} \leq 4 \, \|\ell\|_{L_\alpha^1(\mathbb{R}, \mu_\alpha(dx))} \, \|h\|_{L_\alpha^1(\mathbb{R}, \mu_\alpha(dx))}.
$$
\end{theorem}
\begin{proof}
	According to Statement~(ii) of Theorem~4.4 in~\cite{Bouze2024}, if \( g, h \in L_\alpha^1(\mathbb{R}, \mu_\alpha(dx)) \), then convolution \( g \ast_\alpha h \) is well-defined and belongs to \( L_\alpha^1(\mathbb{R}, \mu_\alpha(dx)) \). Furthermore, the following identity holds
	\begin{equation}
	\label{eq3.1}
	\mathscr{H}_\alpha[g \ast_\alpha h] = \mathscr{H}_\alpha[g] \mathscr{H}_\alpha[h],
	\end{equation}
	Applying $\mathscr{H}_\alpha$-transform on both sides of \eqref{eq3.4}, we obtain $\H_\alpha[f+ (f \ast_\alpha g)]=\mathscr{H}_\alpha[g \ast_\alpha h]$.
	By utilizing identity~\eqref{eq3.1} and the linearity of $\H_\alpha$-transform, this leads to:
	$$
	\mathscr{H}_\alpha[f](\lambda)  \left(1 + \mathscr{H}_\alpha[g](\lambda)\right) = \mathscr{H}_\alpha[g](\lambda)  \mathscr{H}_\alpha[h](\lambda),
	$$
	which means that $$\mathscr{H}_\alpha[f](\lambda) = \frac{\mathscr{H}_\alpha[g](\lambda)}{1+ \mathscr{H}_\alpha[g](\lambda)} \mathscr{H}_\alpha[h](\lambda).$$ Under the condition $1 + \mathscr{H}_\alpha[g](\lambda) \neq 0$ for all $\lambda \in \mathbb{R},$ the above expression is valid. Recall that $L_\alpha^1(\mathbb{R},\mu_\alpha(dx))$ is a Banach algebra under the multiplication convolution \eqref{eq1.7}, (see Remark \ref{daisobanach}). By the Wiener–Lévy invertibility criterion \cite{Naimark},  there exists a function $\ell \in L_\alpha^1(\mathbb{R},\mu_\alpha(dx))$  such that  $$\mathscr{H}_\alpha[\ell](\lambda)=\frac{\mathscr{H}_\alpha[g](\lambda)}{1+ \mathscr{H}_\alpha[g](\lambda)}.$$ This yields $$
	\mathscr{H}_\alpha[f](\lambda) = \mathscr{H}_\alpha[\ell](\lambda) \mathscr{H}_\alpha[h](\lambda) = \mathscr{H}\left[\ell \ast_\alpha h\right](\lambda),
	$$
	which implies, by the injectivity of $\H_\alpha$-transform, then $f(x)=\left(\ell \ast_\alpha h\right)(x)$  almost everywhere on \( \mathbb{R} \). Moreover, since both \( \ell \) and \( h \) belong to \( L_\alpha^1(\mathbb{R}, \mu_\alpha(dx)) \), the convolution \( \ell \ast_\alpha h \) is well-defined and belongs to \( L_\alpha^1(\mathbb{R}, \mu_\alpha(dx)) \). The uniqueness of the solution follows from the fact that the convolution structure \( \ast_\alpha \) is uniquely defined on \( L_\alpha^1 \), as shown in~\cite{Bouze2024}.	
	Finally, by applying Statement~(i) of Theorem~4.4 in~\cite{Bouze2024}, we derive the estimate
	$$
	\|f\|_{L_\alpha^1(\mathbb{R}, \mu_\alpha(dx))} = \|\ell \ast_\alpha h\|_{L_\alpha^1(\mathbb{R}, \mu_\alpha(dx))} \leq 4 \, \|\ell\|_{L_\alpha^1(\mathbb{R}, \mu_\alpha(dx))} \times \|h\|_{L_\alpha^1(\mathbb{R}, \mu_\alpha(dx))}. 
	$$ This completes the proof. 
\end{proof}
\noindent We now examine several special cases concerning the boundedness of the solution that arise from the application of the inequalities established in Section~\ref{sec2}.
\begin{remark} \textup{\textbf{(A).} For $\alpha >0,$ if $f,\ell\in L_\alpha^2(\mathbb{R},\mu_\alpha(dx))$ and $h\in L_\alpha^1(\mathbb{R},\mu_\alpha(dx))$, then $$
		\|f\|_{L_\alpha^2(\mathbb{R},\mu_\alpha(dx))} \leq \sqrt{2}\|\ell\|_{L_\alpha^2(\mathbb{R},\mu_\alpha(dx))} \|h\|_{L_\alpha^1(\mathbb{R},\mu_\alpha(dx))}. 
		$$
\noindent \textbf{(B).}  For $\alpha >0,$ when $\ell, h\in L_\alpha^2(\mathbb{R},\mu_\alpha(dx))$ and $f\in L_\alpha^\infty(\mathbb{R},\mu_\alpha(dx))$, then $$
\|f\|_{L_\alpha^\infty(\mathbb{R},\mu_\alpha(dx))} \leq \sqrt{2}\|\ell\|_{L_\alpha^2(\mathbb{R},\mu_\alpha(dx))} \|h\|_{L_\alpha^2(\mathbb{R},\mu_\alpha(dx))}.
$$           
\noindent \textbf{(C).} Assume that $ p, q, r \in [1, 2] $ satisfy the duality conditions
$$1/p +1/p_1 =1/q +1/q_1 =1/r +1/r_1 =1
\quad
\text{and}\quad 
1/p_1 + 1/q_1 = 1/r_1.
$$
Suppose further that \( f \in L_\alpha^{r_1}(\mathbb{R}, \mu_\alpha(dx)) \), \( \ell \in L_\alpha^p(\mathbb{R}, \mu_\alpha(dx)) \), and \( h \in L_\alpha^q(\mathbb{R}, \mu_\alpha(dx)) \), for \( \alpha > 0 \). Then the solution \( f \) to equation~\eqref{eq3.4} admits the following a priori bound:
$$
\|f\|_{L_\alpha^{r_1}(\mathbb{R}, \mu_\alpha(dx))} \leq ( \sqrt{2})^{\left(\frac{2}{p}+\frac{2}{q}+\frac{2}{r}-3\right)}
\|\ell\|_{L_\alpha^p(\mathbb{R}, \mu_\alpha(dx))}
\|h\|_{L_\alpha^q(\mathbb{R}, \mu_\alpha(dx))}.
$$}	
\end{remark}

\begin{example}[Illustration of Theorem \ref{theo3.1}]
\textup{ Let $\alpha =1/2$, the weighted measure \ref{eq1.5} becomes 
$$\mu_{1 / 2}(d x)=\frac{1}{2^{\frac{1}{2}+1 / 2} \Gamma\left(\frac{1}{2}+1 / 2\right)}|x|^{2 \cdot \frac{1}{2}} d x=\frac{1}{2^1 \Gamma(1)}|x|^1 d x=\frac{1}{2}|x| d x.$$
The kernel defined by \eqref{eq1.3} is expressed as follows
$$\J_{\lambda}(x, 1/2)=\B_0(\lambda x)+\frac{\lambda x}{2} \B_1(\lambda x),$$
where
$
\B_\alpha(x)=\sum_{n=0}^{\infty} \frac{(-1)^n}{n!(\alpha+1)_n}\left(\frac{x}{2}\right)^{2 n} 
$.
When $\alpha=0$ then $\B_0(x)=J_0(x)$ is a Bessel function of the first kind.
When $\alpha=1$ then $\B_1(x)=\frac{2}{x} J_1(x)$, since $J_1(x)=\frac{x}{2} \B_1(x)$. }

\noindent\textup{ This implies that $
\frac{\lambda x}{2} \B_1(\lambda x)=\frac{\lambda x}{2} \cdot \frac{2}{\lambda x} J_1(\lambda x)=J_1(\lambda x).
$ Therefore 
$$\J_{\lambda}(x,1/2)= J_0(\lambda x) + J_1(\lambda x).$$
Define $
g(x)=h(x):=\frac{2 e^{-|x|}}{|x|}. 
$ A direct computation shows
$$
\int_{-\infty}^{\infty}|g(x)| \mu_{1 / 2}(d x)=\int_{-\infty}^{\infty} \frac{2 e^{-|x|}}{|x|} \cdot \frac{1}{2}|x| d x=\int_{-\infty}^{\infty} e^{-|x|} d x=2<\infty,
$$
hence the functions $g,h \in L^1_{1/2} (\R, \mu_{1 / 2}dx)$.\\ By definition \eqref{eq1.6}, we have Hartley-Bessel transform of $g$ is
$$\begin{aligned} 
\mathscr{H}_{1 / 2}[g](\lambda)&=\int_{-\infty}^{\infty} g(x)\J_{\lambda}(x, 1/2) \mu_{1 / 2}(d x)
\\&=\int_{-\infty}^{\infty} e^{-|x|}\big[ J_0(\lambda x)+J_1(\lambda x) \big] d x.
\end{aligned}$$
Because $J_1$ is odd function and $e^{-|x|}$ is even, then the integral of $e^{-|x|}J_1$ vanishes on $\R$. Hence
$$
\mathscr{H}_{1 / 2}[g](\lambda)=2 \int_0^{\infty} e^{-x} J_0(\lambda x) d x.
$$ 
Using the standard Laplace integral $\int_0^{\infty} e^{-a x} J_0(b x) d x=\left(a^2+b^2\right)^{-1 / 2}$ for $a>0$, we deduce that
$$
\mathscr{H}_{1 / 2}[g](\lambda)=\frac{2}{\sqrt{1+\lambda^2}}, \quad \lambda \in \mathbb{R}.
$$
This led to  $
1+ \mathscr{H}_{1 / 2}[g](\lambda)=1+ \frac{2}{\sqrt{1+\lambda^2}}>0, \forall \lambda \in \mathbb{R}
$, hence verification of the non‑vanishing condition of Theorem  \ref{theo3.1} is satisfied. Now, we need to show the existence and uniqueness of the solution. Indeed, by Theorem  \ref{theo3.1} there exists a unique function $\ell \in L_{1 / 2}^1(\mathbb{R}, \mu_{1 / 2}dx)$ such that 
$$
\mathscr{H}_{1 / 2}[\ell](\lambda)=\frac{\mathscr{H}_{1 / 2}[g](\lambda)}{1+\mathscr{H}_{1 / 2}[g](\lambda)}=\frac{2}{\sqrt{1+\lambda^2}+2}
$$
Note that Wiener-L\'evy's invertibility \cite{Naimark} guarantees the existence and uniqueness of $\ell$, because $L^1_{1/2}$ is a commutative Banach algebra equipped with multiplication given by convolution $*_{1/2}$ (see Remark \ref{daisobanach}), and moreover the condition $1+\mathscr{H}_{1 / 2}[g](\lambda)$ does not vanish. The solution $f$ of the integral equation
$$
f(x)+\left(f *_{1 / 2} g\right)(x)=\left(g *_{1 / 2} h\right)(x)
$$
is given by $f=\ell *_{1 / 2} h=\ell *_{1 / 2} g$, since $h=g$. From the factorization property \eqref{eq1.9a} we obtain
$$
\mathscr{H}_{1 / 2}[f](\lambda)=\mathscr{H}_{1 / 2}[\ell](\lambda) \mathscr{H}_{1 / 2}[g](\lambda)=\frac{4}{\sqrt{1+\lambda^2}(\sqrt{1+\lambda^2}+2)} .
$$
A simple algebraic manipulation yields
$
\frac{4}{\sqrt{1+\lambda^2}\left(\sqrt{1+\lambda^2}+2\right)}=\mathscr{H}_{1 / 2}[g](\lambda)-\mathscr{H}_{1 / 2}[\ell](\lambda) .
$
Therefore
$$
\mathscr{H}_{1 / 2}[f]=\mathscr{H}_{1 / 2}[g-\ell],
$$
and by the injectivity of transform $\mathscr{H}_{1 / 2}$ on $L^1_{1 / 2}$, we infer that
$
f=g-\ell \ \text { a.e. on } \mathbb{R} .
$
Consequently $\ell=g-f$ belongs to $L_{1 / 2}^1$ as required. And applying the inequality \eqref{eq1.8} with $p=q=r=1$, then 
$$
\|f\|_{L_{1 / 2}^1}=\left\|\ell *_{1 / 2} g\right\|_{L_{1 / 2}^1} \leq 4\|\ell\|_{L_{1 / 2}^1}\|g\|_{L_{1 / 2}^1}.
$$ }
\end{example}

\vskip 0.3cm
\noindent \textbf{Acknowledgements}\\
We are sincerely indebted to Professor V\~u Kim Tu\^a\'n  for the useful discussions and insightful perspectives he shared on Wiener-L\'evy's invertibility in the context of Banach algebras. We would also like to thank Dr. Trinh Tung for suggesting a beautiful example at the end of the paper, which ensures the validity of the results in Section \ref{sec3}.
\vskip 0.3cm
\noindent \textbf{Disclosure statement}\\
\noindent No potential conflict of interest was reported by the author(s).\\

\noindent \textbf{Funding}\\
\noindent This research received no specific grant from any funding agency.\\
\noindent \textbf{ORCID}\\
\noindent \textsc{Trinh Tuan} { \url{https://orcid.org/0000-0002-0376-0238}}

\end{document}